\documentclass[12pt]{amsart}
\usepackage{amsmath, mathtools}
\usepackage{mathrsfs}
\usepackage{bbm}
\usepackage{amsfonts}
\usepackage{latexsym,amssymb,amscd,amsfonts,epsfig,graphicx}
\usepackage[margin=1.45in]{geometry}

\usepackage{amsthm}
\newtheorem{thm}{Theorem}[section]
\newtheorem{lem}[thm]{Lemma}

\newtheorem{coro}[thm]{Corollary}
\newtheorem{obs}[thm]{Observation}

\def\pf{\noindent{\it Proof.\;\;\:}}
\def\qed{\nopagebreak\hfill{\rule{4pt}{7pt}}
\medbreak}

\DeclareMathOperator{\diam}{diam}
\DeclareMathOperator{\dege}{d}

\theoremstyle{definition}
\newtheorem{defi}[thm]{Definition}

\title[Hadwiger's conjecture for $\ell$-link graphs]{Hadwiger's conjecture for $\ell$-link graphs}

\author[Bin Jia]{Bin Jia}
\address{Department of Mathematics and Statistics\\
The University of Melbourne\\
Melbourne\\
Australia.} \email{jiabinqq@gmail.com}

\thanks{Bin Jia gratefully acknowledges scholarships provided by The University of Melbourne.}

\author[David R. Wood]{David R. Wood}
\address{School of Mathematical Sciences\\
Monash University\\
Melbourne\\
Australia.} \email{david.wood@monash.edu}
\thanks{Research of David Wood is supported by the Australian Research Council}

\begin{document}

\maketitle

\noindent {\bf Abstract.}
In this paper we define and study a new family of graphs that generalises the notions of line graphs and path graphs. Let $G$ be a graph with no loops but possibly with parallel edges. An \emph{$\ell$-link} of $G$ is a walk of $G$ of length $\ell \geqslant 0$ in which consecutive edges are different. We identify an $\ell$-link with its reverse sequence. The \emph{$\ell$-link graph $\mathbb{L}_\ell(G)$} of $G$ is the graph with vertices the $\ell$-links of $G$, such that two vertices are joined by $\mu \geqslant 0$ edges in $\mathbb{L}_\ell(G)$ if they correspond to two subsequences of each of $\mu$ $(\ell + 1)$-links of $G$.

By revealing a recursive structure, we bound from above the chromatic number of $\ell$-link graphs. As a corollary, for a given graph $G$ and large enough $\ell$, $\mathbb{L}_\ell(G)$ is $3$-colourable. By investigating the shunting of $\ell$-links in $G$, we show that the Hadwiger number of a nonempty $\mathbb{L}_\ell(G)$ is greater or equal to that of $G$. Hadwiger's conjecture states that the Hadwiger number of a graph is at least the chromatic number of that graph. The conjecture has been proved by Reed and Seymour (2004) for line graphs, and hence $1$-link graphs. We prove the conjecture for a wide class of $\ell$-link graphs.

\vskip 5pt


\section{Introduction and main results}\label{sec_intro}
We introduce a new family of graphs, called \emph{$\ell$-link graphs}, which generalises the notions of line graphs and path graphs. Such a graph is constructed from a certain kind of walk of length $\ell \geqslant 0$ in a given graph $G$. To ensure that the constructed graph is undirected, $G$ is  undirected, and we identify a walk with its reverse sequence. To avoid loops, $G$ is loopless, and the consecutive edges in each walk are different. Such a walk is called an \emph{$\ell$-link}. For example, a $0$-link is a vertex, a $1$-link is an edge, and a $2$-link consists of two edges with an end vertex in common. An \emph{$\ell$-path} is an $\ell$-link without repeated vertices. We use $\mathscr{L}_\ell(G)$ and $\mathscr{P}_\ell(G)$ to denote the sets of $\ell$-links and $\ell$-paths of $G$ respectively. There have been a number of families of graphs constructed from $\ell$-links. As one of the most commonly studied graphs, the \emph{line graph $\mathbb{L}(G)$}, introduced by Whitney \cite{Whitney1932}, is the simple graph with vertex set $E(G)$, in which two vertices are adjacent if their corresponding edges are incident to a common vertex. More generally, the \emph{$\ell$-path graph $\mathbb{P}_{\ell}(G)$} is the simple graph with vertex set $\mathscr{P}_{\ell}(G)$, where two vertices are adjacent if the union of their corresponding $\ell$-paths forms a path or a cycle of length $\ell + 1$. Note that $\mathbb{P}_{\ell}(G)$ is the $\mathbb{P}_{\ell + 1}$-graph of $G$ introduced by Broersma and Hoede \cite{BH1989}. Inspired by these graphs, we define the \emph{$\ell$-link graph $\mathbb{L}_\ell(G)$} of $G$ to be the graph with vertex set $\mathscr{L}_\ell(G)$, in which two vertices are joined by $\mu \geqslant 0$ edges in $\mathbb{L}_\ell(G)$ if they correspond to two subsequences of each of $\mu$ $(\ell + 1)$-links of $G$. More strict definitions can be found in Section \ref{sec_defiTermi}, together with some other related graphs.

This paper studies the structure, colouring and minors of $\ell$-link graphs including a proof of Hadwiger's conjecture for a wide class of $\ell$-link graphs. By default $\ell \geqslant 0$ is an integer. And all graphs are finite, undirected and loopless. Parallel edges are admitted unless we specify the graph to be \emph{simple}.

\subsection{Graph colouring} Let $t \geqslant 0$ be an integer. A \emph{$t$-colouring} of $G$ is a map $\lambda: V(G) \rightarrow [t] := \{1, 2, \ldots, t\}$ such that $\lambda(u) \neq \lambda(v)$ whenever $u, v \in V(G)$ are adjacent in $G$. A graph with a $t$-colouring is \emph{$t$-colourable}. The \emph{chromatic number $\chi(G)$} is the minimum $t$ such that $G$ is $t$-colourable. Similarly, an \emph{$t$-edge-colouring} of $G$ is a map $\lambda: E(G) \rightarrow [t]$ such that $\lambda(e) \neq \lambda(f)$ whenever $e, f \in E(G)$ are incident to a common vertex in $G$. The \emph{edge-chromatic number $\chi'(G)$} of $G$ is the minimum $t$ such that $G$ admits a $t$-edge-colouring. Let $\chi_\ell(G) := \chi(\mathbb{L}_\ell(G))$, and \emph{$\Delta(G)$} be the maximum degree of $G$. By \cite[Proposition 5.2.2]{Diestel2010}, $\chi_0(G) = \chi(G) \leqslant \Delta(G) + 1$. Shannon \cite{Shannon1949} proved that $\chi_1(G) = \chi'(G) \leqslant \frac{3}{2}\Delta(G)$. We prove a recursive structure for $\ell$-link graphs which leads to the following upper bounds for $\chi_\ell(G)$:

\begin{thm}\label{thm_chrl}
Let $G$ be a graph, $\chi := \chi(G)$, $\chi' := \chi'(G)$, and $\Delta := \Delta(G)$.
\begin{itemize}
\item[\bf (1)] If $\ell \geqslant 0$ is even, then $\chi_\ell(G) \leqslant \min\{\chi, \lfloor(\frac{2}{3})^{\ell/2}(\chi - 3)\rfloor + 3\}$.
\item[\bf (2)] If $\ell \geqslant 1$ is odd, then $\chi_\ell(G) \leqslant \min\{\chi', \lfloor(\frac{2}{3})^{\frac{\ell - 1}{2}}(\chi' - 3)\rfloor + 3\}$.
\item[\bf (3)] If $\ell \neq 1$, then $\chi_\ell(G) \leqslant \Delta + 1$. 
\item[\bf (4)] If $\ell \geqslant 2$, then $\chi_\ell(G) \leqslant \chi_{\ell - 2}(G)$. 
\end{itemize}
\end{thm}

Theorem \ref{thm_chrl} implies that $\mathbb{L}_\ell(G)$ is $3$-colourable for large enough $\ell$.

\begin{coro}\label{coro_chrl} For each graph $G$,
$\mathbb{L}_\ell(G)$ is $3$-colourable in the following cases:
\begin{itemize}
\item[{\bf (1)}] $\ell \geqslant 0$ is even, and either $\chi(G) \leqslant 3$ or $\ell > 2\log_{1.5}(\chi(G) - 3)$.
\item[{\bf (2)}] $\ell \geqslant 1$ is odd, and either $\chi'(G) \leqslant 3$ or $\ell > 2\log_{1.5}(\chi'(G) - 3) + 1$.
\end{itemize}
\end{coro}

As explained in Section \ref{sec_defiTermi}, this corollary is related to and implies a result by Kawai and Shibata \cite{KawaiShibata2002}.

\subsection{Graph minors} By \emph{contracting} an edge we mean identifying its end vertices and deleting possible resulting loops. A graph $H$ is a \emph{minor} of $G$ if $H$ can be obtained from a subgraph of $G$ by contracting edges. An \emph{$H$-minor} is a minor of $G$ that is isomorphic to $H$. The \emph{Hadwiger number $\eta(G)$} of $G$ is the maximum integer $t$ such that $G$ contains a $K_t$-minor. Denote by \emph{$\delta(G)$} the minimum degree of $G$. The \emph{degeneracy $\dege(G)$} of $G$ is the maximum $\delta(H)$ over the subgraphs $H$ of $G$. We prove the following:
\begin{thm}\label{thm_hgeqg}
Let $\ell \geqslant 1$, and $G$ be a graph such that $\mathbb{L}_\ell(G)$ contains at least one edge. Then $\eta(\mathbb{L}_\ell(G)) \geqslant \max\{\eta(G), \dege(G)\}$.
\end{thm}

By definition $\mathbb{L}(G)$ is the underlying simple graph of $\mathbb{L}_1(G)$. And $\mathbb{L}_\ell(G) = \mathbb{P}_\ell(G)$ if $girth(G) > \{\ell, 2\}$.
Thus Theorem \ref{thm_hgeqg} can be applied to path graphs.
\begin{coro}\label{coro_hgeqg}
Let $\ell \geqslant 1$, and $G$ be a graph of girth at least $\ell + 1$ such that $\mathbb{P}_\ell(G)$ contains at least one edge. Then
$\eta(\mathbb{P}_\ell(G)) \geqslant \max\{\eta(G), \dege(G)\}$.
\end{coro}

As a far-reaching generalisation of the four-colour theorem, in 1943, Hugo Hadwiger \cite{Hadwiger1943} conjectured the following:

\medskip
\noindent {\bf Hadwiger's conjecture:} $\eta(G) \geqslant \chi(G)$ for every graph $G$.
\medskip

Hadwiger's conjecture was proved by Robertson, Seymour and Thomas  \cite{RST1993} for $\chi(G) \leqslant 6$. The conjecture for line graphs, or equivalently for $1$-link graphs, was proved by Reed and Seymour \cite{Reed2004}. We prove the following:
\begin{thm}\label{thm_hcllink}
Hadwiger's conjecture is true for $\mathbb{L}_\ell(G)$ in the following cases:
\begin{itemize}
\item[{\bf (1)}] $\ell \geqslant 1$ and $G$ is biconnected. 
\item[{\bf (2)}] $\ell \geqslant 2$ is an even integer. 
\item[{\bf (3)}] $\dege(G) \geqslant 3$ and $\ell > 2\log_{1.5}\frac{\Delta(G) - 2}{\dege(G) - 2} + 3$. 
\item[{\bf (4)}] $\Delta(G) \geqslant 3$ and $\ell > 2\log_{1.5}(\Delta(G) - 2) - 3.83$. 
\item[{\bf (5)}] $\Delta(G) \leqslant 5$.
\end{itemize}
\end{thm}

The corresponding results for path graphs are listed below:
\begin{coro}
Let $G$ be a graph of girth at least $\ell + 1$. Then Hadwiger's conjecture holds for $\mathbb{P}_{\ell}(G)$ in the cases of Theorem \ref{thm_hcllink} (1) -- (5).
\end{coro}

\section{Definitions and terminology}\label{sec_defiTermi}
We now give some formal definitions. A graph $G$ is \emph{null} if $V(G) = \emptyset$, and \emph{nonnull} otherwise. A nonnull graph $G$ is \emph{empty} if $E(G) = \emptyset$, and nonempty otherwise. A \emph{unit} is a vertex or an edge. The subgraph of $G$ induced by $V \subseteq V(G)$ is the maximal subgraph of $G$ with vertex set $V$. And in this case, the subgraph is called an \emph{induced subgraph} of $G$. For $\emptyset \neq E \subseteq E(G)$, the subgraph of $G$ induced by $E \cup V$ is the minimal subgraph of $G$ with edge set $E$, and vertex set including $V$.

For more accurate analysis, we need to define \emph{$\ell$-arcs}. An \emph{$\ell$-arc} (or \emph{$*$-arc} if we ignore the length) of $G$ is an alternating sequence $\vec{L} := (v_0, e_1, \ldots, e_{\ell}, v_\ell)$ of units of $G$ such that the end vertices of $e_i \in E(G)$ are $v_{i - 1}$ and $v_i$ for $i \in [\ell]$, and that $e_i \neq e_{i + 1}$ for $i \in [\ell - 1]$. The \emph{direction} of $\vec{L}$ is its vertex sequence \emph{$(v_0, v_1, \ldots, v_\ell)$}. In algebraic graph theory, $\ell$-arcs in simple graphs have been widely studied \cite{Tutte1947,Tutte1959,Weiss1981,Biggs1993}. Note that $\vec{L}$ and its \emph{reverse $-\vec{L} := (v_\ell, e_{\ell}, \ldots, e_1, v_0)$} are different unless $\ell = 0$. The $\ell$-link (or \emph{$*$-link} if the length is ignored) $L := [v_0, e_1, \ldots, e_{\ell}, v_\ell]$ is obtained by taking $\vec{L}$ and $-\vec{L}$ as a single object. For $0 \leqslant i \leqslant j \leqslant \ell$, the $(j - i)$-arc \emph{$\vec{L}(i, j) := (v_i, e_{i + 1},\ldots, e_{j}, v_j)$} and the $(j - i)$-link \emph{$\vec{L}[i, j] := [v_i, e_{i + 1}, \ldots, e_{j}, v_j]$} are called \emph{segments} of $\vec{L}$ and $L$ respectively. We may write $\vec{L}(j, i) := -\vec{L}(i, j)$, and $\vec{L}[j, i] := \vec{L}[i, j]$. These segments are called \emph{middle segments} if $i + j = \ell$. $L$ is called an \emph{$\ell$-cycle} if $\ell \geqslant 2$, $v_0 = v_\ell$ and $\vec{L}[0, \ell - 1]$ is an $(\ell - 1)$-path. Denote by $\vec{\mathscr{L}}_{\ell}(G)$ and $\mathscr{C}_\ell(G)$ the sets of $\ell$-arcs and $\ell$-cycles of $G$ respectively. Usually, $\vec{e}_i := (v_{i - 1}, e_i, v_i)$ is called an \emph{arc} for short. In particular, $v_0$, $v_\ell$, $e_1$, $e_\ell$, $\vec{e}_1$ and $\vec{e}_\ell$ are called the \emph{tail vertex}, \emph{head vertex}, \emph{tail edge}, \emph{head edge}, \emph{tail arc}, and \emph{head arc} of $\vec{L}$ respectively.

\begin{figure}[!h]
\centering
\includegraphics{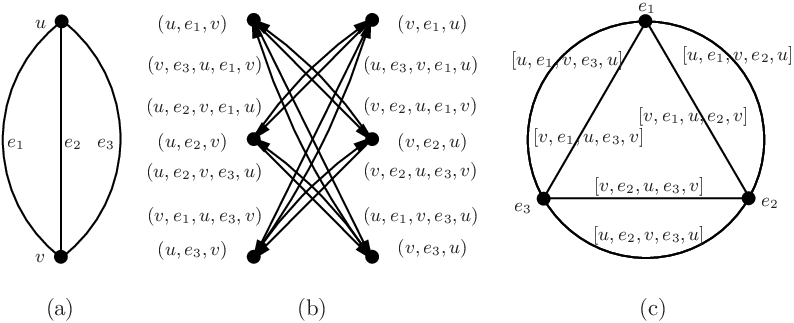}
\caption{(a) $D_3$ $\quad$ (b) $\mathbb{A}_1(D_3)$ $\quad$ (c) $\mathbb{L}_1(D_3)$}
\label{F:D3A1L1}
\end{figure}

Godsil and Royle \cite{GR2001} defined the \emph{$\ell$-arc graph $\mathbb{A}_{\ell}(G)$} to be the digraph with vertex set $\vec{\mathscr{L}}_\ell(G)$, such that there is an arc, labeled by $\vec{Q}$, from $\vec{Q}(0, \ell)$ to $\vec{Q}(1, \ell + 1)$ in $\mathbb{A}_{\ell}(G)$ for every $\vec{Q} \in \vec{L}_{\ell + 1}(G)$. The \emph{$t$-dipole graph $D_t$} is the graph consists of two vertices and $t \geqslant 1$ edges between them. (See Figure \ref{F:D3A1L1}(a) for $D_3$, and Figure \ref{F:D3A1L1}(b) the $1$-arc graph of $D_3$.)
The \emph{$\ell^{th}$ iterated line digraph $\mathbb{A}^\ell(G)$} is $\mathbb{A}_1(G)$ if $\ell = 1$, and $\mathbb{A}_1(\mathbb{A}^{\ell - 1}(G))$ if $\ell \geqslant 2$ (see  \cite{BW1978}). Examples of undirected graphs constructed from $\ell$-arcs can be found in \cite{JLW2010,Jia2011}.

\begin{figure}[!h]
\centering
\includegraphics{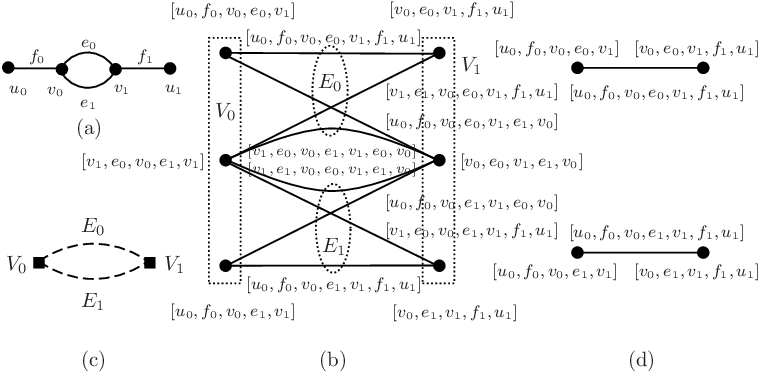}
\caption{(a) $G$ $\quad$ (b) $H := \mathbb{L}_2(G)$ $\quad$ (c) $H_{(\mathcal{V}, \mathcal{E})}$ $\quad$ (d) $\mathbb{P}_2(G)$}
\label{F:GL2P2Quotient}
\end{figure}

\emph{Shunting} of $\ell$-arcs was introduced by Tutte \cite{Tutte1966}. We extend this motion to $\ell$-links. For $\ell, s \geqslant 0$, and $\vec{Q} \in \vec{\mathscr{L}}_{\ell + s}(G)$, let $\vec{L}_i := \vec{Q}(i, \ell + i)$ for $i \in [0, s]$, and $\vec{Q}_i := \vec{L}(i - 1, \ell + i)$ for $i \in [s]$. Let \emph{$Q^{[\ell]} := [L_0, Q_1, L_1, \ldots, L_{s - 1}, Q_{s}, L_{s}]$}. We say $L_0$ can be \emph{shunted} to $L_s$ through $\vec{Q}$ or $Q$. \emph{$Q^{\{\ell\}} := \{L_0, L_1, \ldots, L_s\}$} is the set of \emph{images} during this shunting. For $L, R \in \mathscr{L}_\ell(G)$, we say $L$ can be \emph{shunted} to $R$ if there are $\ell$-links $L = L_0, L_1, \ldots, L_s = R$ such that $L_{i - 1}$ can be shunted to $L_{i}$ through some $*$-arc $\vec{Q}_i$ for $i \in [s]$. In Figure \ref{F:GL2P2Quotient}, $[u_0, f_0, v_0, e_0, v_1]$ can be shunted to $[v_1, e_0, v_0, e_1, v_1]$ through $(u_0, f_0, v_0, e_0, v_1, f_1, u_1)$ and $(u_1, f_1, v_1, e_0, v_0, e_1, v_1)$.

For $L, R \in \mathscr{L}_\ell(G)$ and $\mathscr{Q} \subseteq \mathscr{L}_{\ell + 1}(G)$, denote by \emph{$\mathscr{Q}(L, R)$} the set of $Q \in \mathscr{Q}$ such that $L$ can be shunted to $R$ through $Q$. We show in Section \ref{sec_genStruc} that $|\mathscr{Q}(L, R)|$ is $0$ or $1$ if $G$ is simple, and can be up to $2$ if $\ell \geqslant 1$ and $G$ contains parallel edges. A more formal definition of $\ell$-link graphs is given below:
\begin{defi}\label{defi_llinkgraph}
Let $\mathscr{L} \subseteq \mathscr{L}_\ell(G)$, and $\mathscr {Q} \subseteq \mathscr{L}_{\ell + 1}(G)$. The \emph{partial $\ell$-link graph $\mathbb{L}(G, \mathscr{L}, \mathscr{Q})$}  of $G$, with respect to $\mathscr{L}$ and $\mathscr{Q}$, is the graph with vertex set $\mathscr{L}$, such that $L, R \in \mathscr{L}$ are joined by exactly $|\mathscr{Q}(L, R)|$ edges. In particular, $\mathbb{L}_\ell(G) = \mathbb{L}(G, \mathscr{L}_\ell(G), \mathscr{L}_{\ell + 1}(G))$ is the $\ell$-link graph of $G$.
\end{defi}
\noindent {\bf Remark.} We assign exclusively to each edge of $\mathbb{L}_\ell(G)$ between $L, R \in \mathscr{L}_\ell(G)$ a $Q \in \mathscr{L}_{\ell + 1}(G)$ such that $L$ can be shunted to $R$ through $Q$, and refer to this edge simply as $Q$. In this sense, $Q^{[\ell]} := [L, Q, R]$ is a $1$-link of $\mathbb{L}_{\ell}(G)$.

For example, the $1$-link graph of $D_3$ can be seen in Figure \ref{F:D3A1L1}(c). A $2$-link graph is given in Figure \ref{F:GL2P2Quotient}(b), and a $2$-path graph is depicted in Figure \ref{F:GL2P2Quotient}(d).

Reed and Seymour \cite{Reed2004} pointed out that proving Hadwiger's conjecture for line graphs of multigraphs is more difficult than for that of simple graphs. This motivates us to work on the $\ell$-link graphs of multigraphs. Diestel \cite[page 28]{Diestel2010} explained that, in some situations, it is more natural to develop graph theory for multigraphs. The observation below follows from the definitions:

\begin{obs}
$\mathbb{L}_0(G) = G$, $\mathbb{P}_1(G) = \mathbb{L}(G)$, and $\mathbb{P}_{\ell}(G)$ is the underlying simple graph of $\mathbb{L}_\ell(G)$ for $\ell \in \{0, 1\}$. For $\ell \geqslant 2$, $\mathbb{P}_{\ell}(G) = \mathbb{L}(G, \mathscr{P}_{\ell}(G), \mathscr{P}_{\ell + 1}(G)\\ \cup \mathscr{C}_{\ell + 1}(G))$ is an induced subgraph of $\mathbb{L}_{\ell}(G)$. If $G$ is simple, then $\mathbb{P}_{\ell}(G) = \mathbb{L}_\ell(G)$ for $\ell \in \{0, 1, 2\}$. Further, $\mathbb{P}_{\ell}(G) = \mathbb{L}_\ell(G)$ if $girth(G) > \max\{\ell, 2\}$.
\end{obs}

Let $\vec{Q} \in \vec{\mathscr{L}}_{\ell + s}(G)$, and $[L_0, Q_1, L_1, \ldots, L_{s - 1}, Q_{s}, L_{s}] := Q^{[\ell]}$. From Definition \ref{defi_llinkgraph}, for $i \in [s]$, $Q_i$ is an edge of $H := \mathbb{L}_\ell(G)$ between $L_{i - 1}, L_i \in V(H)$. So $Q^{[\ell]}$ is an $s$-link of $H$. In Figure \ref{F:GL2P2Quotient}(b), $[u_0, f_0, v_0, e_0, v_1, e_1, v_0, e_0, v_1]^{[2]} = [[u_0, f_0, v_0, e_0, v_1], [u_0, f_0, v_0, e_0, v_1, e_1, v_0], [v_0, e_0, v_1, e_1, v_0], [v_0, e_0, v_1, e_1, v_0, e_0, v_1]$,\\ $[v_1, e_1, v_0, e_0, v_1]]$ is a $2$-path of $H$.

We say $H$ is \emph{homomorphic} to $G$, written $H \rightarrow G$, if there is an injection $\alpha: V(H) \cup E(H) \rightarrow V(G) \cup E(G)$ such that for $w \in V(H)$, $f \in E(H)$ and $[u, e, v] \in \mathscr{L}_1(H)$, their images $w^\alpha \in V(G)$, $f^\alpha \in E(G)$ and $[u^\alpha, e^\alpha, v^\alpha] \in \mathscr{L}_1(G)$. In this case,  $\alpha$ is called a \emph{homomorphism} from $H$ to $G$. The definition here is a generalisation of the one for simple graphs by Godsil and Royle \cite[Page 6]{GR2001}. A bijective homomorphism is an \emph{isomorphism}. By Hell and Ne{\v{s}}et{\v{r}}il \cite{HellNesetril2004}, $\chi(H) \leqslant \chi(G)$ if $H \rightarrow G$. For instance, $\vec{L} \mapsto L$ for $\vec{L} \in \vec{\mathscr{L}}_{\ell}(G) \cup \vec{\mathscr{L}}_{\ell + 1}(G)$ can be seen as a homomorphism from $\mathbb{A}_\ell(G)$ to $\mathbb{L}_\ell(G)$. By Bang-Jensen and Gutin \cite{BangGutin2009}, $\mathbb{A}_\ell(G) \cong \mathbb{A}^\ell(G)$. So $\chi(\mathbb{A}^\ell(G)) = \chi(\mathbb{A}_\ell(G)) \leqslant \chi(\mathbb{L}_\ell(G)) \leqslant \chi_\ell(G)$. We emphasize that $\chi(\mathbb{A}^\ell(G))$ might be much less than $\chi_\ell(G)$. For example, as depicted in Figure \ref{F:D3A1L1}, when $t \geqslant 3$, $\chi(\mathbb{A}^\ell(D_t)) = 2 < t = \chi_\ell(D_t)$. Kawai and Shibata proved that $\mathbb{A}^\ell(G)$ is $3$-colourable for large enough $\ell$. By the analysis above, Corollary \ref{coro_chrl} implies this result.

A graph homomorphism from $H$ is usually represented by a vertex partition $\mathcal {V}$ and an edge partition $\mathcal {E}$ of $H$ such that: {\bf (a)} each part of $\mathcal {V}$ is an independent set of $H$, and {\bf (b)} each part of $\mathcal {E}$ is incident to exactly two parts of $\mathcal {V}$. In this situation, for different $U, V \in \mathcal {V}$, define \emph{$\mu(U, V)$} to be the number of parts of $\mathcal {E}$ incident to both $U$ and $V$. The \emph{quotient graph $H_{(\mathcal {V}, \mathcal {E})}$} of $H$ is defined to be the graph with vertex set $\mathcal {V}$, and for every pair of different $U, V \in \mathcal {V}$, there are exactly $\mu(U, V)$ edges between them. To avoid ambiguity, for $V \in \mathcal {V}$ and $E \in \mathcal {E}$, we use $V_{\mathcal {V}}$ and $E_{\mathcal {E}}$ to denote the corresponding vertex and edge of $H_{(\mathcal {V}, \mathcal {E})}$, which defines a graph homomorphism from $H$ to $H_{(\mathcal {V}, \mathcal {E})}$. Sometimes, we only need the underlying simple graph \emph{$H_{\mathcal {V}}$} of $H_{(\mathcal {V}, \mathcal {E})}$.

For $\ell \geqslant 2$, there is a natural partition in an $\ell$-link graph. For each $R \in \mathscr{L}_{\ell - 2}(G)$, let $\mathscr{L}_\ell(R)$ be the set of $\ell$-links of $G$ with middle segment $R$. Clearly, $\mathcal {V}_{\ell}(G) := \{\mathscr{L}_{\ell}(R) \neq \emptyset | R \in \mathscr{L}_{\ell - 2}(G)\}$ is a vertex partition of $\mathbb{L}_\ell(G)$. And $\mathcal {E}_\ell(G) := \{\mathscr{L}_{\ell + 1}(P) \neq \emptyset | P \in \mathscr{L}_{\ell - 1}(G)\}$ is an edge partition of $\mathbb{L}_\ell(G)$. Consider the $2$-link graph $H$ in Figure \ref{F:GL2P2Quotient}(b). The vertex and edge partitions of $H$ are indicated by the dotted rectangles and ellipses respectively. The corresponding quotient graph is given in Figure \ref{F:GL2P2Quotient}(c).

Special partitions are required to describe the structure of $\ell$-link graphs. Let $H$ be a graph admitting partitions $\mathcal {V}$ of $V(H)$ and $\mathcal {E}$ of $E(H)$ that satisfy {\bf (a)} and {\bf (b)} above. $(\mathcal {V}, \mathcal {E})$ is called an \emph{almost standard partition} of $H$ if further:

{\bf (c)} each part of $\mathcal {E}$ induces a complete bipartite subgraph of $H$,

{\bf (d)} each vertex of $H$ is incident to at most two parts of $\mathcal {E}$,

{\bf (e)} for each $V \in \mathcal{V}$, and different $E, F \in \mathcal {E}$, $V$ contains at most one vertex incident to both $E$ and $F$.

If $\ell \geqslant 2$ is an even integer, and $G$ is a simple graph, then $\mathbb{L}_\ell(G)$ is isomorphic to the \textit{$(2, \ell/2)$-double star graph} of $G$ introduced by Jia \cite{Jia2011}. While this paper focuses on the combinatorial properties including connectedness, colouring and minors of $\mathbb{L}_\ell(G)$, a series of companion papers have been composed to contribute to the recognition and determination problems and algorithms. For example, a joint work by Ellingham and Jia \cite{EllinghamJia2016} shows that, for a given graph $H$, there is at most one pair $(G, \ell)$, where $\ell \geqslant 2$, and $G$ is a simple graph of minimum degree at least $3$, such that $\mathbb{L}_\ell(G)$ is isomorphic to $H$. Moreover, such a pair can be determined from $H$ in linear time.
\section{General structure of $\ell$-link graphs}\label{sec_genStruc}
We begin by determining some basic properties of $\ell$-link graphs, including their multiplicity and connectedness. The work in this section forms the basis for our main results on colouring and minors of $\ell$-link graphs.

Let us first fix some concepts by two observations.

\begin{obs}\label{obs_llinkregu}
The number of edges of $\mathbb{L}_{\ell}(G)$ is equal to the number of vertices of $\mathbb{L}_{\ell + 1}(G)$. In particular, if $G$ is $r$-regular for some $r \geqslant 2$, then this number is $|E(G)|(r - 1)^{\ell}$. If further $\ell \geqslant 1$, then $\mathbb{L}_{\ell}(G)$ is $2(r - 1)$-regular.
\end{obs}
\pf Let $G$ be $r$-regular, $n := |V(G)|$ and $m := |E(G)|$. We prove that $|\mathscr{L}_{\ell + 1}(G)| = m(r - 1)^{\ell}$ by induction on $\ell$. It is trivial for $\ell = 0$. For $\ell = 1$, $|\mathscr{L}_2(v)| = \binom{r}{2}$, and hence $|\mathscr{L}_{2}(G)| = \binom{r}{2}n = m(r - 1)$. Inductively assume $|\mathscr{L}_{\ell - 1}(G)| = m(r - 1)^{\ell - 2}$ for some $\ell \geqslant 2$. For each $R \in \mathscr{L}_{\ell - 1}(G)$, we have $|\mathscr{L}_{\ell + 1}(R)| = (r - 1)^2$ since $r \geqslant 2$. Thus $|\mathscr{L}_{\ell + 1}(G)| = |\mathscr{L}_{\ell - 1}(G)|(r - 1)^2 = m(r - 1)^{\ell}$ as desired. The other assertions follow from the definitions.
\qed

\begin{obs}\label{obs_llinkbipar}
Let $n, m \geqslant 2$. If $\ell \geqslant 1$ is odd, then $\mathbb{L}_\ell(K_{n, m})$ is $(n + m - 2)$-regular with order $nm[(n - 1)(m - 1)]^{\frac{\ell - 1}{2}}$. If $\ell \geqslant 2$ is even, then $\mathbb{L}_\ell(K_{n, m})$ has average degree $\frac{4(n - 1)(m - 1)}{n + m - 2}$, and order $\frac{1}{2}nm(n + m - 2)[(n - 1)(m - 1)]^{\frac{\ell}{2} - 1}$.
\end{obs}
\pf
Let $\ell \geqslant 1$ be odd, and $L$ be an $\ell$-link of $K_{n, m}$ with middle edge incident to a vertex $u$ of degree $n$ in $K_{n, m}$. It is not difficult to see that $L$ can be shunted in one step to $n - 1$ $\ell$-links whose middle edge is incident to $u$. By symmetry, each vertex of $\mathbb{L}_\ell(K_{n, m})$ is incident to $(n - 1) + (m - 1) = n + m - 2$ edges. Now we prove $|\mathscr{L}_\ell(K_{n, m})| = nm[(n - 1)(m - 1)]^{\frac{\ell - 1}{2}}$ by induction on $\ell$. Clearly, $|\mathscr{L}_1(K_{n, m})| = |E(K_{n, m})| = nm$. Inductively assume $|\mathscr{L}_{\ell - 2}(K_{n, m})| = nm[(n - 1)(m - 1)]^{\frac{\ell - 3}{2}}$ for some $\ell \geqslant 3$. For each $R \in \mathscr{L}_{\ell - 2}(K_{n, m})$, we have $|\mathscr{L}_{\ell}(R)| = (n - 1)(m - 1)$. So $|\mathscr{L}_\ell(K_{n, m})| = |\mathscr{L}_{\ell - 2}(K_{n, m})|(n - 1)(m - 1) = nm[(n - 1)(m - 1)]^{\frac{\ell - 1}{2}}$ as desired. The even $\ell$ case is similar.
\qed

\subsection{Loops and multiplicity}
Our next observation is a prerequisite for the study of the chromatic number since it indicates that $\ell$-link graphs are loopless.

\begin{obs}\label{obs_loopless}
For each $(\ell + 1)$-arc $\vec{Q}$, we have $\vec{Q}[0, \ell] \neq \vec{Q}[1, \ell + 1]$.
\end{obs}
\pf Let $G$ be a graph, and $\vec{Q} := (v_0, e_1, \ldots, e_{\ell + 1}, v_{\ell + 1}) \in \vec{\mathscr{L}}_{\ell + 1}(G)$. Since $G$ is loopless, $v_0 \neq v_1$ and hence $\vec{Q}(0, \ell) \neq \vec{Q}(1, \ell + 1)$. So the statement holds for $\ell = 0$. Now let $\ell \geqslant 1$. Suppose for a contradiction that $\vec{Q}(0, \ell) = -\vec{Q}(1, \ell + 1)$. Then $v_i = v_{\ell + 1 - i}$ and $e_{i + 1} = e_{\ell + 1 - i}$ for $i \in \{0, 1, \ldots, \ell\}$. If $\ell = 2s$ for some integer $s \geqslant 1$, then $v_s = v_{s + 1}$, contradicting that $G$ is loopless.  If $\ell = 2s + 1$ for some $s \geqslant 0$, then $e_{s + 1} = e_{s + 2}$, contradicting the definition of a $*$-arc.
\qed

The following statement indicates that, for each $\ell \geqslant 1$, $\mathbb{L}_\ell(G)$ is simple if $G$ is simple, and has multiplicity exactly $2$ otherwise.
\begin{obs}\label{obs_mul2}
Let $G$ be a graph, $\ell \geqslant 1$, and $L_0, L_1 \in \mathscr{L}_{\ell}(G)$. Then $L_0$ can be shunted to $L_1$ through two $(\ell + 1)$-links of $G$ if and only if $G$ contains a $2$-cycle $O := [v_0, e_0, v_1, e_1, v_0]$, such that one of the following cases holds:
\begin{itemize}
\item[{\bf (1)}] $\ell \geqslant 1$ is odd, and $L_i = [v_i, e_i, v_{1 - i}, e_{1 - i}, \ldots, v_i, e_i, v_{1 - i}] \in \mathscr{L}_\ell(O)$ for $i \in \{0, 1\}$. In this case, $[v_i, e_i, v_{1 - i}, e_{1 - i}, \ldots, v_{1 - i}, e_{1 - i}, v_i] \in \mathscr{L}_{\ell + 1}(O)$, for $i \in \{0, 1\}$, are the only two $(\ell + 1)$-links available for the shunting.

\item[{\bf (2)}] $\ell \geqslant 2$ is even, and $L_i = [v_i, e_i, v_{1 - i}, e_{1 - i}, \ldots, v_{1 - i}, e_{1 - i}, v_i] \in \mathscr{L}_\ell(O)$ for $i \in \{0, 1\}$. In this case, $[v_i, e_i, v_{1 - i}, e_{1 - i}, \ldots, v_i, e_i, v_{1 - i}] \in \mathscr{L}_{\ell + 1}(O)$, for $i \in \{0, 1\}$, are the only two $(\ell + 1)$-links available for the shunting.
\end{itemize}
\end{obs}
\pf
$(\Leftarrow)$ is trivial. For $(\Rightarrow)$, since $L_0$ can be shunted to $L_1$, there exists $\vec{L} := (v_0, e_0, v_1, \ldots, v_\ell, e_\ell, v_{\ell + 1}) \in \vec{\mathscr{L}}_{\ell + 1}(G)$ such that $L_i = \vec{L}[i, \ell + i]$ for $i \in \{0, 1\}$. Let $\vec{R} \in \vec{\mathscr{L}}_{\ell + 1}(G) \setminus \{\vec{L}\}$ such that $L_i = \vec{R}[i, \ell + i]$. Then $\vec{L}(i, \ell + i)$ equals $\vec{R}(i, \ell + i)$ or $\vec{R}(\ell + i, i)$. Suppose for a contradiction that $\vec{L}(0, \ell) = \vec{R}(0, \ell)$. Then $\vec{L}(1, \ell) = \vec{R}(1, \ell)$. Since $\vec{L} \neq \vec{R}$, we have $\vec{L}(1, \ell + 1) \neq \vec{R}(1, \ell + 1)$. Thus $\vec{L}(1, \ell + 1) = \vec{R}(\ell + 1, 1)$, and hence $\vec{L}(2, \ell + 1) = \vec{R}(\ell, 1) = \vec{L}(\ell, 1)$, contradicting Observation \ref{obs_loopless}. So $\vec{L}(0, \ell) = \vec{R}(\ell, 0)$. Similarly, $\vec{L}(1, \ell) = \vec{R}(\ell + 1, 1)$. Consequently, $\vec{L}(0, \ell - 1) = \vec{R}(\ell, 1) = \vec{L}(2, \ell + 1)$; that is, $v_j = v_0$ and $e_j = e_0$ if $j \in [0, \ell]$ is even, while $v_j = v_1$ and $e_j = e_1$ if $j \in [0, \ell + 1]$ is odd.
\qed

\subsection{Connectedness} This subsection characterises when $\mathbb{L}_\ell(G)$ is connected.
A middle segment of $L \in \mathscr{L}_\ell(G)$ is a \emph{middle unit}, written \emph{$c_L$}, if it is a unit of $G$. Note that $c_L$ is a vertex if $\ell$ is even, and is an edge otherwise. Denote by \emph{$G(\ell)$} the subgraph of $G$ induced by the middle units of $\ell$-links of $G$.

The lemma below is important in dealing with the connectedness of $\ell$-link graphs. Before stating it, we define a \emph{conjunction} operation, which is an extension of an operation by Biggs \cite[Chapter 17]{Biggs1993}. Let $\vec{L} := (v_0, e_1, v_1, \ldots, e_\ell, v_\ell) \in \vec{\mathscr{L}}_\ell(G)$ and $\vec{R} := (u_0, f_1, u_1, \ldots, f_s, u_s) \in \vec{\mathscr{L}}_s(G)$ such that $v_\ell = u_0$ and $e_{\ell} \neq f_1$. The \emph{conjunction} of  $\vec{L}$ and $\vec{R}$ is $(\vec{L}.\vec{R}) := (v_0, e_1, \ldots, e_{\ell}, v_{\ell} = u_0, f_1, \ldots, f_s, u_s) \in \vec{\mathscr{L}}_{\ell + s}(G)$ or $[\vec{L}.\vec{R}] := [v_0, e_1, \ldots, e_{\ell}, v_{\ell} = u_0, f_1, \ldots, f_s, u_s] \in \mathscr{L}_{\ell + s}(G)$.

\begin{lem}\label{lem_GlConnect}
Let $\ell, s \geqslant 0$, and $G$ be a connected graph. Then $G(\ell)$ is connected. And each $s$-link of $G(\ell)$ is a middle segment of a $(2 \lfloor\frac{\ell}{2}\rfloor + s)$-link of $G$. Moreover, for $\ell$-links $L$ and $R$ of $G$, there is an $\ell$-link $L'$ with middle unit $c_L$, and an $\ell$-link $R'$ with middle unit $c_R$, such that $L'$ can be shunted to $R'$.
\end{lem}
\pf For $\ell \in \{0, 1\}$, since $G$ is connected, $G(\ell) = G$ and the lemma holds. Let $\ell := 2m \geqslant 2$ be even. $u, v \in V(G(\ell))$ if and only if they are middle vertices of some $\vec{L}, \vec{R} \in \vec{\mathscr {L}}_\ell(G)$ respectively. Since $G$ is connected, there exists some $\vec{P} \in \vec{\mathscr{L}}_s(G)$ from $(u, e, u_1)$ to $(v_{s - 1}, f, v)$. By Observation \ref{obs_loopless}, $\vec{L}[m - 1, m] \neq \vec{L}[m, m + 1]$. For such an $s$-arc $\vec{P}$, without loss of generality, $e \neq \vec{L}[m - 1, m]$, and similarly, $f \neq \vec{R}[m, m + 1]$. Then $\vec{P}$ is a middle segment of $\vec{Q} := (\vec{L}(0, m).\vec{P}.\vec{R}(m, 2m)) \in \vec{\mathscr{L}}_{\ell + s}(G)$. So $\vec{P} \in \vec{\mathscr{L}}_{s}(G(\ell))$. And $L' := \vec{Q}[0, \ell]$ can be shunted to $R' := \vec{Q}[s, \ell + s]$ through $\vec{Q}$. The odd $\ell$ case is similar.
\qed

Sufficient conditions for $\mathbb{A}_\ell(G)$ to be strongly connected
can be found in \cite[Page 76]{GR2001}. The following corollary of Lemma \ref{lem_GlConnect} reveals a strong relationship between the shunting of $\ell$-links and the connectedness of $\ell$-link graphs.
\begin{coro}\label{coro_LlConnect}
For a connected graph $G$, $\mathbb{L}_\ell(G)$ is connected if and only if any two $\ell$-links of $G$ with the same middle unit can be shunted to each other.
\end{coro}

We now present our main result of this section, which plays a key role in dealing with the graph minors of $\ell$-link graphs in Section \ref{sec_Hadwiger}.
\begin{lem}\label{lem_hub}
Let $G$ be a graph, and $X$ be a connected subgraph of $G(\ell)$. Then for every pair of $\ell$-links $L$ and $R$ of $X$, $L$ can be shunted to $R$ under the restriction that in each step, the middle unit of the image of $L$ belongs to $X$.
\end{lem}
\pf
First we consider the case that $c_L$ is in $R$. Then there is a common segment $Q$ of $L$ and $R$ of maximum length containing $c_L$. Without loss of generality, assign directions to $L$ and $R$ such that $\vec{L} = (\vec{L}_0.\vec{Q}.\vec{L}_1)$ and $\vec{R} = (\vec{R}_1.\vec{Q}. \vec{R}_0)$, where $\vec{L}_i \in \vec{\mathscr{L}}_{\ell_i}(X)$ and $\vec{R}_i \in \vec{\mathscr{L}}_{s_i}(X)$ for $i \in \{0, 1\}$ such that $s_1 \geqslant s_0$. Then $\ell \geqslant \ell_0 + \ell_1 = s_0 + s_1 \geqslant s_1$. Let $x$ be the head vertex and $e$  be the head edge of $\vec{L}$. Since $c_L$ is in $Q$, $\ell_0 \leqslant \ell/2$. Since $X$ is a subgraph of $G(\ell)$, by Lemma \ref{lem_GlConnect}, there exists $\vec{L}_2 \in \vec{\mathscr{L}}_{\ell_0}(G)$ with tail vertex $x$ and tail edge different from $e$. Let $y$ be the tail vertex and $f$ be the tail edge of $\vec{R}$. Then there exits $\vec{R}_2 \in \vec{\mathscr{L}}_{s_0}(G)$ with head vertex $y$ and head edge different from $f$. We can shunt $L$ to $R$ first through $(\vec{L}.\vec{L}_2) \in \vec{\mathscr{L}}_{\ell + \ell_0}(G)$, then $-(\vec{R}_2.\vec{R}_1.\vec{Q}.\vec{L}_1.\vec{L}_2) \in \vec{\mathscr{L}}_{\ell + \ell_0 + \ell_1}(G)$, and finally $(\vec{R}_2.\vec{R}) \in \vec{\mathscr{L}}_{\ell + s_0}(G)$. Since $\ell_0 \leqslant \ell/2$ and $s_0 \leqslant s_1 \leqslant \ell/2$, the middle unit of each image is inside $L$ or $R$.

Secondly, we consider the case that $c_L$ is not in $R$. Then there exists a segment $Q$ of $L$ of maximum length that contains $c_L$, and is edge-disjoint with $R$. Since $X$ is connected, there exists a shortest $*$-arc $\vec{P}$ from a vertex $v$ of $R$ to a vertex $u$ of $L$. Then $P$ is edge-disjoint with $Q$ because of its minimality. Without loss of generality, assign directions to $L$ and $R$ such that $u$ separates $\vec{L}$ into $(\vec{L}_0.\vec{L}_1)$ with $c_L$ on $L_1$, and $v$ separates $\vec{R}$ into $(\vec{R}_1.\vec{R}_0)$, where $L_i$ is of length $\ell_i$ while $R_i$ is of length $s_i$ for $i \in \{0, 1\}$, such that $s_1 \geqslant s_0$. Then $\ell_0, s_0 \leqslant \ell/2$. Let $x$ be the head vertex and $e$ be the head edge of $\vec{L}$. Since $\ell_0 \leqslant \ell/2$ and $X$ is a subgraph of $G(\ell)$, by Lemma \ref{lem_GlConnect}, there exists an $\ell_0$-arc $\vec{L}_2$ of $G$ with tail vertex $x$ and tail edge different from $e$. Let $y$ be the tail vertex and $f$ be the tail edge of $\vec{R}$. Then there exits an $s_0$-arc $\vec{R}_2$ of $G$ with head vertex $y$ and head edge different from $f$. Now we can shunt $L$ to $R$ through $(\vec{L}.\vec{L}_2)$, $-(\vec{R}_2.\vec{R}_1.\vec{P}.\vec{L}_1.\vec{L}_2)$ and $(\vec{R}_2.\vec{R})$ consecutively. One can check that in this process the middle unit of each image belongs to $L, P$ or $R$.
\qed

From Lemma \ref{lem_hub}, the set of $\ell$-links of a connected $G(\ell)$ serves as a `hub' in the shunting of $\ell$-links of $G$. More explicitly, for $L, R \in \mathscr{L}_\ell(G)$, if we can shunt $L$ to $L' \in \mathscr{L}_\ell(G(\ell))$, and $R$ to $R' \in \mathscr{L}_\ell(G(\ell))$, then $L$ can be shunted to $R$ since $L'$ can be shunted to $R'$. Thus we have the following corollary which provides a more efficient way to test the connectedness of $\ell$-link graphs.
\begin{coro}\label{coro_hub} Let $G$ be a graph. Then
$\mathbb{L}_\ell(G)$ is connected if and only if $G(\ell)$ is connected, and each $\ell$-link of $G$ can be shunted to an $\ell$-link of $G(\ell)$.
\end{coro}

\section{Chromatic number of $\ell$-link graphs}\label{sec_chromatic}
In this section, we reveal a recursive structure of $\ell$-link graphs, which leads to an upper bound for the chromatic number of $\ell$-link graphs.

\begin{lem}\label{lem_VEpartition}
Let $G$ be a graph and $\ell \geqslant 2$ be an integer. Then
$(\mathcal {V}, \mathcal {E}) := (\mathcal {V}_{\ell}(G), \mathcal {E}_{\ell}(G))$ is an almost standard partition of $H := \mathbb{L}_\ell(G)$. Further, $H_{(\mathcal {V}, \mathcal {E})}$ is isomorphic to an induced subgraph of $\mathbb{L}_{\ell - 2}(G)$.
\end{lem}
\pf First we verify that $(\mathcal {V}, \mathcal {E})$ is an almost standard partition of $H$.

{\bf (a)} We prove that, for each $R \in \mathscr{L}_{\ell - 2}(G)$, $V := \mathscr{L}_\ell(R) \in \mathcal{V}$ is an independent set of $H$. Suppose not. Then there are $\vec{L}, \vec{L}' \in \vec{\mathscr{L}}_\ell(G)$ such that $L, L' \in V$, and $L$ can be shunted to $L'$ in one step. Then $R = \vec{L}[1, \ell - 1]$ can be shunted to $R = \vec{L}'[1, \ell - 1]$ in one step, contradicting Observation \ref{obs_loopless}.

{\bf (b)} Here we show that each $E \in \mathcal {E}$ is incident to exactly two parts of $\mathcal {V}$. By definition there exists $P \in \mathscr {L}_{\ell - 1}(G)$ with $\mathscr{L}_{\ell + 1}(P) = E$. Let $\{L, R\} := P^{\{\ell - 2\}}$. Then $\mathscr{L}_{\ell}(L)$ and $\mathscr{L}_{\ell}(R)$ are the only two parts of $\mathcal {V}$ incident to $E$.

{\bf (c)} We explain that each $E \in \mathcal{E}$ is the edge set of a complete bipartite subgraph of $H$. By definition there exists $\vec{P} \in \vec{\mathscr {L}}_{\ell - 1}(G)$ with $\mathscr{L}_{\ell + 1}(P) = E$. Let $A := \{[\vec{e}.\vec{P}] \in \mathscr{L}_\ell(G)\}$ and $B := \{[\vec{P}.\vec{f}] \in \mathscr{L}_\ell(G)\}$. One can check that $E$ induces a complete bipartite subgraph of $H$ with bipartition $A \cup B$.

{\bf (d)} We prove that each $v \in V(H)$ is incident to at most two parts of $\mathcal{E}$. By definition there exists $Q \in \mathscr{L}_\ell(G)$ with $Q = v$. Then the set of edge parts of $\mathcal{E}$ incident to $v$ is $\{\mathscr {L}_{\ell + 1}(L) \neq \emptyset| L \in Q^{\{\ell - 1\}}\}$ with cardinality at most $2$.

{\bf (e)} Let $v$ be a vertex of $V \in \mathcal{V}$ incident to different $E, F \in \mathcal {E}$. We explain that $v$ is uniquely determined by $V$, $E$ and $F$.
By definition there exists $\vec{P} \in \vec{\mathscr {L}}_{\ell - 2}(G)$ such that $V = \mathscr {L}_{\ell}(P)$. There also exists $Q := [\vec{e}_1.\vec{P}.\vec{e}_{\ell}] \in \mathscr {L}_{\ell}(P)$ such that $v = Q$. Besides, there are $L, R \in \mathscr{L}_{\ell - 1}(G)$ such that $E = \mathscr{L}_{\ell + 1}(L)$ and $F = \mathscr{L}_{\ell + 1}(R)$. Then $\{L, R\} = Q^{\{\ell - 1\}}$ since $L \neq R$. Note that $Q$ is uniquely determined by $Q^{\{\ell - 1\}}$ and $c_Q = c_P$. Thus it is uniquely determined by $E = \mathscr{L}_{\ell + 1}(L), F = \mathscr{L}_{\ell + 1}(R)$ and $V = \mathscr {L}_{\ell}(P)$.

Now we show that $H_{(\mathcal {V}, \mathcal {E})}$ is isomorphic to an induced subgraph of $\mathbb{L}_{\ell - 2}(G)$. Let $X$ be the subgraph of $\mathbb{L}_{\ell - 2}(G)$ of vertices $L \in \mathscr{L}_{\ell - 2}(G)$ such that
$\mathscr {L}_{\ell}(L) \neq \emptyset$, and edges $Q \in \mathscr {L}_{\ell - 1}(G)$ such that $\mathscr {L}_{\ell + 1}(Q) \neq \emptyset$. One can check that $X$ is an induced subgraph of $\mathbb{L}_{\ell - 2}(G)$. An isomorphism from $H_{(\mathcal {V}, \mathcal {E})}$ to $X$ can be defined as the injection sending $\mathscr {L}_{\ell}(L) \neq \emptyset$ to $L$, and $\mathscr {L}_{\ell + 1}(Q) \neq \emptyset$ to $Q$.
\qed

Below we give an interesting algorithm for colouring a class of graphs.
\begin{lem}\label{lem_colourr}
Let $H$ be a graph with a $t$-colouring such that each vertex of $H$ is adjacent to at most $r \geqslant 0$ differently coloured vertices. Then $\chi(H) \leqslant \lfloor \frac{tr}{r + 1} \rfloor + 1$.
\end{lem}
\pf
The result is trivial for $t = 0$ since, in this case, $\chi(H) = 0$. If $r + 1 \geqslant t \geqslant 1$, then $\lfloor \frac{tr}{r + 1} \rfloor + 1 = t$, and the lemma holds since $t \geqslant \chi(H)$.

Now assume $t \geqslant r + 2 \geqslant 2$. Let $U_1, U_2, \ldots, U_t$ be the colour classes of the given colouring. For $i \in [t]$, denote by $i$ the colour assigned to vertices in $U_i$. Run the following algorithm: For $j = 1, \ldots, t$, and for each $u \in U_{t - j + 1}$, let $s \in [t]$ be the minimum integer that is not the colour of a neighbour of $u$ in $H$; if $s < t - j + 1$, then recolour $u$ by $s$.

In the algorithm above, denote by $C_i$ the set of colours used by the vertices in $U_i$ for $i \in [t]$. Let $k := \lfloor\frac{t - 1}{r + 1}\rfloor$. Then $t - 1 \geqslant k(r + 1) \geqslant k \geqslant 1$. We claim that after $j \in [0, k]$ steps, $C_{t - i + 1} \subseteq [ir + 1]$ for $i \in [j]$, and $C_i = \{i\}$ for $i \in [t - j]$. This is trivial for $j = 0$. Inductively assume it holds for some $j \in [0, k - 1]$. In the $(j + 1)^{th}$ step, we change the colour of each $u \in U_{t - j}$ from $t - j$ to the minimum $s \in [t]$ that is not used by the neighbourhood of $u$. It is enough to show that $s \leqslant (j + 1)r + 1$.

First suppose that all neighbours of $u$ are in $\bigcup_{i \in [t - j - 1]}U_i$. By the analysis above, $t - j - 1 \geqslant t - k \geqslant kr + 1 \geqslant r + 1$. So at least one part of $\mathcal{S} := \{U_i| i \in [t - j - 1]\}$ contains no neighbour of $u$. From the induction hypothesis, $C_i = \{i\}$ for $i \in [t - j - 1]$. Hence at least one colour in $[r + 1]$ is not used by the neighbourhood of $u$; that is, $s \leqslant r + 1 \leqslant (j + 1)r + 1$.

Now suppose that $u$ has at least one neighbour in $\bigcup_{i \in [t - j + 1, t]}U_i$. By the induction hypothesis, $\bigcup_{i \in [t - j + 1, t]}C_i \subseteq [jr + 1]$. At the same time, $u$ has neighbours in at most $r - 1$ parts of $\mathcal{S}$. So the colours possessed by the neighbourhood of $u$ are contained in $[jr + 1 + r - 1] = [(j + 1)r]$. Thus $s \leqslant (j + 1)r + 1$. This proves our claim.

The claim above indicates that, after the $k^{th}$ step, $C_{t - i + 1} \subseteq [ir + 1]$ for $i \in [k]$, and $C_i = \{i\}$ for $i \in [t - k]$. Hence we have a $(t - k)$-colouring of $H$ since $t - k \geqslant kr + 1$. Therefore, $\chi(H) \leqslant t - k = \lceil\frac{tr + 1}{r + 1}\rceil = \lfloor \frac{tr}{r + 1} \rfloor + 1$.
\qed

Lemma \ref{lem_VEpartition} indicates that $\mathbb{L}_\ell(G)$ is homomorphic to $\mathbb{L}_{\ell - 2}(G)$ for $\ell \geqslant 2$. So by \cite[Proposition 1.1]{Cameron2006}, $\chi_\ell(G) \leqslant \chi_{\ell - 2}(G)$. By Lemma \ref{lem_VEpartition}, every vertex of $\mathbb{L}_\ell(G)$ has neighbours in at most two parts of $\mathcal {V}_\ell(G)$, which enables us to improve the upper bound on $\chi_\ell(G)$.

\begin{lem}\label{lem_colourr}
Let $G$ be a graph, and $\ell \geqslant 2$. Then $\chi_\ell(G) \leqslant \lfloor \frac{2}{3}\chi_{\ell - 2}(G)\rfloor + 1$.
\end{lem}
\pf
By Lemma \ref{lem_VEpartition}, $(\mathcal {V}, \mathcal {E}) := (\mathcal {V}_{\ell}(G), \mathcal {E}_{\ell}(G))$ is an almost standard partition of $H := \mathbb{L}_\ell(G)$. So each vertex of $H$ has neighbours in at most two parts of $\mathcal{V}$. Further, $H_{\mathcal{V}}$ is a subgraph of $\mathbb{L}_{\ell - 2}(G)$. So $\chi_\ell(G) \leqslant \chi := \chi(H_{\mathcal{V}}) \leqslant \chi_{\ell - 2}(G)$.

We now construct a $\chi$-colouring of $H$ such that each vertex of $H$ is adjacent to at most two differently coloured vertices. By definition $H_{\mathcal{V}}$ admits a $\chi$-colouring with colour classes $K_1, \ldots, K_{\chi}$. For $i \in [\chi]$, assign the colour $i$ to each vertex of $H$ in $U_i := \bigcup_{V_{\mathcal {V}} \in K_i}V$. One can check that this is a desired colouring. In Lemma \ref{lem_colourr}, letting $t = \chi$ and $r = 2$ yields that $\chi_\ell(G) \leqslant \lfloor\frac{2}{3}\chi\rfloor + 1$. Recall that $\chi \leqslant \chi_{\ell - 2}(G)$. Thus the lemma follows. \qed

As shown below, Lemma \ref{lem_colourr} can be applied recursively to produce an upper bound for $\chi_\ell(G)$ in terms of $\chi(G)$ or $\chi'(G)$.

\medskip
\noindent{\textbf{Proof of Theorem \ref{thm_chrl}}}. When $\ell \in \{0, 1\}$, it is trivial for {\bf (1)}{\bf (2)} and {\bf (4)}.
By \cite[Proposition 5.2.2]{Diestel2010}, $\chi_0 = \chi \leqslant \Delta + 1$. So {\bf (3)} holds. Now let $\ell \geqslant 2$. By Lemma \ref{lem_VEpartition}, $H := \mathbb{L}_\ell(G)$ admits an almost standard partition $(\mathcal{V}, \mathcal{E}) := (\mathcal{V}_\ell(G), \mathcal{E}_\ell(G))$, such that $H_{(\mathcal{V}, \mathcal{E})}$ is an induced subgraph of $\mathbb{L}_{\ell - 2}(G)$. By definition each part of $\mathcal{V}$ is an independent set of $H$. So $H \rightarrow \mathbb{L}_{\ell - 2}(G)$, and $\chi_\ell \leqslant \chi_{\ell - 2}$. This proves {\bf (4)}. Moreover, each vertex of $H$ has neighbours in at most two parts of $\mathcal{V}$. By Lemma \ref{lem_colourr}, $\chi_{\ell} := \chi_\ell(G) \leqslant \frac{2\chi_{\ell - 2}}{3} + 1$. Continue the analysis, we have
$\chi_\ell \leqslant \chi_{\ell - 2i}$, and $\chi_\ell - 3 \leqslant (\frac{2}{3})^{i}(\chi_{\ell - 2i} - 3)$ for $1 \leqslant i \leqslant \lfloor \ell/2\rfloor$. Therefore, if $\ell$ is even, then $\chi_\ell \leqslant \chi_0 = \chi \leqslant \Delta + 1$, and $\chi_\ell - 3 \leqslant  (\frac{2}{3})^{\ell/2}(\chi - 3)$. Thus {\bf (1)} holds. Now let $\ell \geqslant 3$ be odd. Then $\chi_\ell \leqslant \chi_1 = \chi'$, and $\chi_\ell - 3 \leqslant (\frac{2}{3})^{\frac{\ell - 1}{2}}(\chi' - 3)$. This verifies {\bf (2)}. As a consequence, $\chi_\ell \leqslant \chi_3 \leqslant \frac{2}{3}(\chi' - 3) + 3 = \frac{2}{3}\chi' + 1$. By Shannon \cite{Shannon1949}, $\chi' \leqslant \frac{3}{2}\Delta$. So $\chi_\ell \leqslant \Delta + 1$, and hence {\bf (3)} holds.
\qed

The following corollary of Theorem \ref{thm_chrl} implies that Hadwiger's conjecture is true for $\mathbb{L}_\ell(G)$ if $G$ is regular and $\ell \geqslant 4$.
\begin{coro}\label{coro_chrHadl}
Let $G$ be a graph with $\Delta := \Delta(G) \geqslant 3$. Then $\chi_\ell(G) \leqslant 3$ for all $\ell > 2\log_{1.5}(\Delta - 2) + 3$. Further, Hadwiger's conjecture holds for $\mathbb{L}_\ell(G)$ if $\ell > 2\log_{1.5}(\Delta - 2) - 3.83$, or $\dege := \dege(G) \geqslant 3$ and $\ell > 2\log_{1.5}\frac{\Delta - 2}{\dege - 2} + 3$.
\end{coro}
\pf
By Theorem \ref{thm_chrl}, for each $t \geqslant 3$, $\chi_\ell := \chi_\ell(G) \leqslant t$ if $(\frac{2}{3})^{\ell/2}(\Delta - 2)  < t - 2$ and $(\frac{2}{3})^{\frac{\ell - 1}{2}}(\frac{3}{2}\Delta - 3) < t - 2$. Solving these inequalities gives $\ell > 2\log_{1.5}(\Delta - 2) - 2\log_{1.5}(t - 2) + 3$. Thus $\chi_\ell \leqslant 3$ if $\ell > 2\log_{1.5}(\Delta - 2) + 3$. So the first statement holds. By Robertson et al. \cite{RST1993} and Theorem \ref{thm_hgeqg}, Hadwiger's conjecture holds for $\mathbb{L}_\ell(G)$ if $\ell \geqslant 1$ and $\chi_\ell \leqslant \max\{6, \dege\}$. Letting $t = 6$ gives that $\ell > 2\log_{1.5}(\Delta - 2) - 4\log_{1.5}2 + 3$. Letting $t = \dege \geqslant 3$ gives that $\ell > 2\log_{1.5}\frac{\Delta - 2}{\dege - 2} + 3$. So the corollary holds since $4\log_{1.5}2 - 3 > 3.83$.
\qed

\noindent{\textbf{Proof of Theorem \ref{thm_hcllink}(3)(4)(5)}}. {\bf (3)} and {\bf (4)} follow from Corollary \ref{coro_chrHadl}. Now consider {\bf (5)}. By Reed and Seymour \cite{Reed2004}, Hadwiger's conjecture holds for $\mathbb{L}_1(G)$. If $\ell \geqslant 2$ and $\Delta \leqslant 5$, by Theorem \ref{thm_chrl}(3), $\chi_\ell(G) \leqslant 6$. In this case, Hadwiger's conjecture holds for $\mathbb{L}_\ell(G)$ by Robertson et al.  \cite{RST1993}.\qed

\section{Complete minors of $\ell$-link graphs}\label{sec_Hadwiger}
It has been proved in the last section that Hadwiger's conjecture is true for $\mathbb{L}_{\ell}(G)$ if $\ell$ is large enough. In this section, we further investigate the minors, especially the complete minors, of $\ell$-link graphs. To see the intuition of our method, let $v$ be a vertex of degree $t$ in $G$. Then $\mathbb{L}_1(G)$ contains a $K_t$-subgraph whose vertices correspond to the edges of $G$ incident to $v$. For $\ell \geqslant 2$, roughly speaking, we extend $v$ to a subgraph $X$ of diameter less than $\ell$, and extend each edge incident to $v$ to an $\ell$-link of $G$ starting from a vertex of $X$. By studying the shunting of these $\ell$-links, we find a $K_t$-minor in $\mathbb{L}_\ell(G)$.

For subgraphs $X, Y$ of $G$, let \emph{$\vec{E}(X, Y)$} be the set of arcs of $G$ from $V(X)$ to $V(Y)$, and \emph{$E(X, Y)$} be the set of edges of $G$ between $V(X)$ and $V(Y)$.

\begin{lem}\label{lem_ltConnected}
Let $\ell \geqslant 1$ be an integer, $G$ be a graph, and $X$ be a subgraph of $G$ with $\diam(X) < \ell$ such that $Y := G - V(X)$ is connected. If $t := |E(X, Y)| \geqslant 2$, then $\mathbb{L}_\ell(G)$ contains a $K_t$-minor.
\end{lem}
\pf Let $\vec{e}_1, \ldots, \vec{e}_t$ be distinct arcs in $\vec{E}(Y, X)$. Say $\vec{e}_i = (y_i, e_i, x_i)$ for $i \in [t]$. Since $\diam(X) < \ell$, there is a dipath  $\vec{P}_{ij}$ of $X$ from $x_i$ to $x_j$ of length $\ell_{ij} \leqslant \ell - 1$ such that $P_{ij} = P_{ji}$. Since $Y$ is connected, it contains a  dipath $\vec{Q}_{ij}$ from $y_i$ to $y_j$. Since $t \geqslant 2$, $O_i := [\vec{P}_{i\; i'}.-\vec{e}_{i'}.\vec{Q}_{i'\; i}. \vec{e}_i]$ is a cycle of $G$, where $i' := (i \mod t) + 1$. Thus $H := \mathbb{L}_\ell(G)$ contains a cycle $\mathbb{L}_\ell(O_1)$, and hence a $K_2$-minor. Now let $t \geqslant 3$, and $\vec{L}_i \in \vec{\mathscr{L}}_\ell(O_i)$ with head arc $\vec{e}_i$. Then $[\vec{L}_i.\vec{P}_{ij}]^{[\ell]} \in \mathscr{L}_{\ell_{ij}}(H)$. And the union of the units of $[\vec{L}_i.\vec{P}_{ij}]^{[\ell]}$ over $j \in [t]$ is a connected subgraph $X_i$ of $H$. In the remainder of the proof, for distinct $i, j \in [t]$, we show that $X_i$ and $X_j$ are disjoint. Further, we construct a path in $H$ between $X_i$ and $X_j$ that is internally disjoint with its counterparts, and has no inner vertex in any of $V(X_1), \ldots, V(X_t)$. Then by contracting each $X_i$ into a vertex, and each path into an edge, we obtain a $K_t$-minor of $H$.

First of all, assume for a contradiction that there are different $i, j \in [t]$ such that $X_i$ and $X_j$ share a common vertex that corresponds to an $\ell$-link $R$ of $G$. Then by definition, there exists some $p \in [t]$ such that $R$ can be obtained by shunting $L_i$ along $(\vec{L}_i.\vec{P}_{ip})$ by some $s_i \leqslant \ell_{ip}$ steps. So $R = [\vec{L}_i(s_i, \ell).\vec{P}_{ip}(0, s_i)]$. Similarly, there are $q \in [t]$ and $s_j \leqslant \ell_{jq}$ such that $R = [\vec{L}_j(s_j, \ell).\vec{P}_{jq}(0, s_j)]$. Recall that $E(X) \cap E(X, Y) = E(Y) \cap E(X, Y) = \emptyset$. So $e_i = \vec{L}_i[\ell - 1, \ell]$ and $e_j = \vec{L}_j[\ell - 1, \ell]$ belong to both $L_i$ and $L_j$. By the definition of $O_i$, this happens if and only if $i = j'$ and $j = i'$, which is impossible since $t \geqslant 3$.

Secondly, for different $i, j \in [t]$, we define a path of $H$ between $X_i$ and $X_j$. Clearly, $L_i$ can be shunted to $L_j$ through $\vec{R}_{ij}' := (\vec{L}_i.\vec{P}_{ij}.-\vec{L}_j)$ in $G$. In this shunting, $L_i' := [\vec{L}_i(\ell_{ij}, \ell).\vec{P}_{ij}]$ is the last image corresponding to a vertex of $X_i$, while $L_j' := [\vec{P}_{ij}.\vec{L}_j(\ell, \ell_{ij})]$ is the first image corresponding to a vertex of $X_j$. Further, $L_i'$ can be shunted to $L_j'$ through $\vec{R}_{ij} := (\vec{L}_i(\ell_{ij}, \ell).\vec{P}_{ij}.\vec{L}_j(\ell, \ell_{ij})) \in \vec{\mathscr{L}}_{2\ell - \ell_{ij}}(G)$, which is a subsequence of $\vec{R}_{ij}'$. Then $R_{ij}^{[\ell]}$ is an $(\ell - \ell_{ij})$-path of $H$ between $X_i$ and $X_j$. We show that for each $p \in [t]$, $X_p$ contains no inner vertex of $R_{ij}^{[\ell]}$. When $\ell - \ell_{ij} = 1$, $R_{ij}^{[\ell]}$ contains no inner vertex. Now assume $\ell - \ell_{ij} \geqslant 2$. Each inner vertex of $R_{ij}^{[\ell]}$ corresponds to some $Q_{ij} := [\vec{L}_i(s_i, \ell).\vec{P}_{ij}.\vec{L}_j(\ell, \ell + \ell_{ij} - s_i)] \in \mathscr{L}_\ell(G)$, where $\ell_{ij} + 1 \leqslant s_i \leqslant \ell - 1$. Assume for a contradiction that for some $p \in [t]$, $X_p$ contains a vertex corresponding to $Q_{ij}$. By definition there exists $q \in [t]$ such that $Q_{ij} = [\vec{L}_p(s_p, \ell).\vec{P}_{pq}(0, s_p)]$, where $0 \leqslant s_p \leqslant \ell_{pq}$. Without loss of generality, $(\vec{L}_i(s_i, \ell).\vec{P}_{ij}.\vec{L}_j(\ell, \ell + \ell_{ij} - s_i)) = (\vec{L}_p(s_p, \ell).\vec{P}_{pq}(0, s_p))$. Since $e_j$ and $e_p$ are not in $P_{pq}$, hence $\vec{e}_j$ belongs to $-\vec{L}_p$ and $\vec{e}_p$ belongs to $-\vec{L}_j$. By the definition of $\vec{L}_i$, this happens only when $j = p'$ and $p = j'$, contradicting $t \geqslant 3$.

We now show that $R_{ij}^{[\ell]}$ and $R_{pq}^{[\ell]}$ are internally disjoint, where $i \neq j$, $p \neq q$ and $\{i, j\} \neq \{p, q\}$. Suppose not. Then by the analysis above, there are $s_i$ and $s_p$ with $\ell_{ij} + 1 \leqslant s_i \leqslant \ell - 1$ and $\ell_{pq} + 1 \leqslant s_p \leqslant \ell - 1$ such that $Q_{ij} = Q_{pq}$. Without loss of generality, $(\vec{L}_i(s_i, \ell).\vec{P}_{ij}.\vec{L}_j(\ell, \ell + \ell_{ij} - s_i)) = (\vec{L}_p(s_p, \ell).\vec{P}_{pq}.\vec{L}_q(\ell, \ell + \ell_{pq} - s_p))$. If $s_i = s_p$, then $\vec{e}_i = \vec{e}_p$ and $\vec{e}_j = \vec{e}_q$ since $E(X) \cap E(X, Y) = \emptyset$; that is, $i = p$ and $j = q$, contradicting $\{i, j\} \neq \{p, q\}$. Otherwise, with no loss of generality, $s_i > s_p$. Then $\vec{e}_q$ and $\vec{e}_i$ belong to $\vec{L}_j$ and $\vec{L}_p$ respectively; that is, $i = p$ and $j = q$, again contradicting $\{i, j\} \neq \{p, q\}$.

In summary, $X_1, \ldots, X_t$ are vertex-disjoint connected subgraphs, which are pairwise connected by internally disjoint  $*$-links $R_{ij}^{[\ell]}$ of $H$, such that no inner vertex of $R_{ij}^{[\ell]}$ is in $V(X_1)\cup \ldots \cup V(X_t)$. So by contracting each $X_i$ to a vertex, and $R_{ij}^{[\ell]}$ to an edge, we obtain a $K_t$-minor of $H$. \qed

\begin{lem}\label{lem_ltConCycle}
Let $\ell \geqslant 1$, $G$ be a graph, and $X$ be a subgraph of $G$ with $\diam(X) < \ell$ such that $Y := G - V(X)$ is connected and contains a cycle. Let $t := |E(X, Y)|$. Then $\mathbb{L}_\ell(G)$ contains a $K_{t + 1}$-minor.
\end{lem}
\pf Let $O$ be a cycle of $Y$. Then $H := \mathbb{L}_\ell(G)$ contains a cycle $\mathbb{L}_\ell(O)$ and hence a $K_2$-minor. Now assume $t \geqslant 2$. Let $\vec{e}_1, \ldots, \vec{e}_t$ be distinct arcs in $\vec{E}(Y, X)$. Say $\vec{e}_i = (y_i, e_i, x_i)$ for $i \in [t]$. Since $Y$ is connected, there is a dipath $\vec{P}_i$ of $Y$ of minimum length $s_i \geqslant 0$ from some vertex $z_i$ of $O$ to $y_i$. Let $\vec{Q}_i$ be an $\ell$-arc of $O$ with head vertex $z_i$. Then $\vec{L}_i := (\vec{Q}_i.\vec{P}_i.\vec{e}_i)(s_i + 1, \ell + s_i + 1) \in \vec{\mathscr{L}}_\ell(G)$. Since $\diam(X) \leqslant \ell - 1$, there is a dipath $\vec{P}_{ij}$ of $X$ of length $\ell_{ij} \leqslant \ell - 1$ from $x_i$ to $x_j$ such that $P_{ij} = P_{ji}$.

Clearly, $[\vec{L}_i.\vec{P}_{ij}]^{[\ell]}$ is an $\ell_{ij}$-link of $H$. And the union of the units of $[\vec{L}_i.\vec{P}_{ij}]^{[\ell]}$ over $j \in [t]$ induces a connected subgraph $X_i$ of $H$. For different $i, j \in [t]$, let $R_{ij} := [\vec{L}_i(\ell_{ij}, \ell).\vec{P}_{ij}.\vec{L}_j(\ell, \ell_{ij})] = R_{ji} \in \mathscr{L}_{2\ell - \ell_{ij}}(G)$. Then $R_{ij}^{[\ell]}$ is an $(\ell - \ell_{ij})$-path of $H$ between $X_i$ and $X_j$. As in the proof of Lemma \ref{lem_ltConnected}, it is easy to check that $X_1, \ldots, X_t$ are vertex-disjoint connected subgraphs of $H$, which are pairwise connected by internally disjoint paths $R_{ij}^{[\ell]}$. Further, no inner vertex of $R_{ij}^{[\ell]}$ is in $V(X_1)\cup \ldots \cup V(X_t)$. So a $K_t$-minor of $H$ is obtained accordingly.

Finally, let $Z$ be the connected subgraph of $H$ induced by the units of $\mathbb{L}_\ell(O)$ and $[\vec{Q}_i.\vec{P}_i]^{[\ell]}$ over $i \in [t]$. Then $Z$ is vertex-disjoint with $X_i$ and with the paths $R_{ij}^{[\ell]}$. Moreover, $Z$ sends an edge $(\vec{Q}_i.\vec{P}_i.\vec{e}_i)(s_i, \ell + s_i + 1)^{[\ell]}$ to each $X_i$. Thus $H$ contains a $K_{t+1}$-minor. \qed

In the following, we use the `hub' (described after Lemma \ref{lem_hub}) to construct certain minors in $\ell$-link graphs.
\begin{coro}\label{coro_lminor}
Let $\ell \geqslant 0$, $G$ be a graph, $M$ be a minor of $G(\ell)$ such that each branch set contains an $\ell$-link. Then $\mathbb{L}_{\ell}(G)$ contains an $M$-minor.
\end{coro}
\pf Let $X_1, \ldots, X_t$ be the branch sets of an $M$-minor of $G(\ell)$ such that $X_i$ contains an $\ell$-link for each $i \in [t]$.
For any connected subgraph $Y$ of $G(\ell)$ contains at least one $\ell$-link, let $\mathbb{L}_{\ell}(G, Y)$ be the subgraph of $H := \mathbb{L}_\ell(G)$ induced by the $\ell$-links of $G$ of which the middle units are in $Y$. Let $H(Y)$ be the union of the components of $\mathbb{L}_\ell(G, Y)$ which contains at least one vertex corresponding to an $\ell$-link of $Y$. By Lemma \ref{lem_hub}, $H(Y)$ is connected.

By definition each edge of $M$ corresponds to an edge $e$ of $G(\ell)$ between two different branch sets, say $X_i$ and $X_j$. Let $Y$ be the graph consisting of $X_i, X_j$ and $e$. Then $H(X_i)$ and $H(X_j)$ are vertex-disjoint since $X_i$ and $X_j$ are vertex-disjoint. By the analysis above, $H(X_i)$ and $H(X_j)$ are connected subgraphs of the connected graph $H(Y)$. Thus there is a path $Q$ of $H(Y)$ joining $H(X_i)$ and $H(X_j)$ only at end vertices. Further, if $\ell$ is even, then $Q$ is an edge; otherwise, $Q$ is a $2$-path whose middle vertex corresponds to an $\ell$-link $L$ of $Y$ such that $c_L = e$. This implies that $Q$ is internally disjoint with its counterparts and has no inner vertex in any branch set. Then, by contracting each $H(X_i)$ to a vertex, and $Q$ to an edge, we obtain an $M$-minor of $H$.
\qed

Now we are ready to give a lower bound for the Hadwiger number of $\mathbb{L}_\ell(G)$.

\medskip
\noindent{\textbf{Proof of Theorem \ref{thm_hgeqg}}}. Since $H := \mathbb{L}_\ell(G)$ contains an edge, $t := \eta(H) \geqslant 2$.
We first show that $t \geqslant \dege := \dege(G)$. By definition there exists a subgraph $X$ of $G$ of $\delta(X) = \dege$. We may assume that $\dege \geqslant 3$. Then $X$ contains an $(\ell - 1)$-link $P$ such that $\mathscr{L}(P) \neq \emptyset$. By Lemma \ref{lem_VEpartition}, $\mathscr{L}^{[\ell]}(P)$ is the edge set of a complete bipartite subgraph of $H$ with a $K_{\dege - 1, \dege - 1}$-subgraph. By Zelinka \cite{Zelinka1976}, $K_{\dege - 1, \dege - 1}$ contains a $K_{\dege}$-minor. Thus $t \geqslant \dege$ as desired.

We now show that
$t \geqslant \eta := \eta(G)$. If $\eta = 3$, then $G$ contains a cycle $O$ of length at least $3$, and $H$ contains a $K_3$-minor contracted from $\mathbb{L}_\ell(O)$. Now assume that $G$ is connected with $\eta \geqslant 4$. Repeatedly delete vertices of degree $1$ in $G$ until $\delta(G) \geqslant 2$. Then $G = G(\ell)$. Clearly, this process does not reduce the Hadwiger number of $G$. So $G$ contains branch sets of a $K_\eta$-minor covering $V(G)$ (see \cite{Wood2011}). If every branch set contains an $\ell$-link, then the statement follows from Corollary \ref{coro_lminor}. Otherwise, there exists some branch set $X$ with $\diam(X) < \ell$. Since $\eta \geqslant 4$, $Y := G - V(X)$ is connected and contains a cycle. Thus by Lemma \ref{lem_ltConCycle}, $H$ contains a $K_\eta$-minor since $|E(X, Y)| \geqslant \eta - 1$.
\qed

Here we prove Hadwiger's conjecture for $\mathbb{L}_\ell(G)$ for even $\ell \geqslant 2$.

\medskip
\noindent{\textbf{Proof of Theorem \ref{thm_hcllink}(2)}.} Let $\dege := \dege(G)$, $\ell \geqslant 2$ be an even integer, and $H := \mathbb{L}_{\ell}(G)$. By \cite[Proposition 5.2.2]{Diestel2010}, $\chi := \chi(G) \leqslant \dege + 1$. So by Theorem \ref{thm_chrl}, $\chi(H) \leqslant \min\{\dege + 1, \frac{2}{3}\dege + \frac{5}{3}\}$. If $\dege \leqslant 4$, then $\chi(H) \leqslant 5$. By Robertson et al. \cite{RST1993}, Hadwiger's conjecture holds for $H$ in this case. Otherwise, $\dege \geqslant 5$. By Theorem \ref{thm_hgeqg}, $\eta(H) \geqslant \dege \geqslant \frac{2}{3}\dege + \frac{5}{3} \geqslant \chi(H)$ and the statement follows.
\qed

We end this paper by proving Hadwiger's conjecture for $\ell$-link graphs of biconnected graphs for $\ell \geqslant 1$.

\medskip
\noindent{\textbf{Proof of Theorem \ref{thm_hcllink}(1)}.} By Reed and Seymour \cite{Reed2004}, Hadwiger's conjecture holds for $H := \mathbb{L}_\ell(G)$ for $\ell = 1$. By Theorem \ref{thm_hcllink}(2), the conjecture is true if $\ell \geqslant 2$ is even. So we only need to consider the situation that $\ell \geqslant 3$ is odd. If $G$ is a cycle, then $H$ is a cycle and the conjecture holds \cite{Hadwiger1943}. Now let $v$ be a vertex of $G$ with degree $\Delta := \Delta(G) \geqslant 3$. By Theorem \ref{thm_chrl}, $\chi(H) \leqslant \Delta + 1$. Since $G$ is biconnected, $Y := G - v$ is connected. By Lemma \ref{lem_ltConCycle}, if $Y$ contains a cycle, then $\eta(H) \geqslant \Delta + 1 \geqslant \chi(H)$. Now assume that $Y$ is a tree, which implies that $G$ is $K_4$-minor free. By Lemma \ref{lem_ltConnected}, $\eta(H) \geqslant \Delta$. By Theorem \ref{thm_chrl}, $\chi(H) \leqslant \chi' := \chi'(G)$. So it is enough to show that $\chi' = \Delta$.

Let $U := \{u \in V(Y)|\;\deg_Y(u) \leqslant 1\}$. Then $|U| \geqslant \Delta(Y)$. Let $\hat{G}$ be the underlying simple graph of $G$, $t := \deg_{\hat{G}}(v) \geqslant 1$ and $\hat{\Delta} := \Delta(\hat{G}) \geqslant t$. Since $G$ is biconnected, $U \subseteq N_G(v)$. So $t \geqslant |U| \geqslant \Delta(Y)$. Let $u \in U$. When $|U| = 1$, $t = \deg_{\hat{G}}(u) = 1$. When $|U| \geqslant 2$, $\deg_{\hat{G}}(u) = 2 \leqslant |U| \leqslant t$. Thus $t = \hat{\Delta}$. Juvan et al. \cite{JMT1999} proved that the edge-chromatic number of a $K_4$-minor free simple graph equals the maximum degree of this graph. So $\hat{\chi}' := \chi'(\hat{G}) = \hat{\Delta}$ since $\hat{G}$ is simple and $K_4$-minor free. Note that all parallel edges of $G$ are incident to $v$. So $\chi' = \hat{\chi}' + \deg_G(v) - t = \hat{\Delta} + \Delta - \hat{\Delta} = \Delta$ as desired.
\qed

\bibliographystyle{plain}
\bibliography{Xbib}
\end{document}